\documentstyle [14pt]{article}
\begin{document}
\title{
Free coherent states
\\
and $p$-adic numbers
}
\author{S.V.Kozyrev}
\maketitle

\footnotetext{
e-mail: kozyrev@class.mian.su,\quad Address: OSV, Inst. Chem. Phys.,117334,
Kosygina 4, Moscow}

\begin{abstract}
Free coherent states for a system with two degrees of freedom is defined.
Existence of the homeomorphism
of the ring of integer 2-adic numbers to the set of coherent states
corresponding to an eigenvalue of the  operator of annihilation is proved.
It is shown that the
metric of free Fock space induces the 2-adic topology on the set of coherent
states.
\end{abstract}

Free (or Boltzmannian)
Fock space has been considered in some recent works
on quantum chromodynamics
\cite{MasterFld}, \cite{GopGross}, \cite{DougLi}
and noncommutative probability
\cite{AccLu}, \cite{Maassen}.

The subject of this work is free coherent states.
We will introduce free coherent states and investigate
the metric space of coherent states corresponding to a fixed eigenvalue
of the operator of annihilation.
It will be shown that
archimedean metric of Hilbert space
induce non-archimedean metric on the set of free coherent states.
The main result of present paper is the construction of
the homeomorphism from the ring of integer 2-adic numbers to
the space of free coherent states with topology defined
by the Hilbert metric. We will consider the system with two degrees of freedom.
The system with one degree of freedom was investigated in \cite{KozAAV}.

The free commutation relations are particular case of
$q$-deformed relations
$$A_i A_j^{\dag}-qA_j^{\dag}A_i=\delta_{ij}  $$
with $q=0$.
A correspondence of  $q$-deformed commutation relations and
non-archimedean (ultrametric) geometry was discussed in
\cite{ArVqGPnAG}.
Non-archimedean mathematical physics
was studied in \cite{VVZ}.

Free coherent states lies in
the free Fock space. Free  (or Boltzmannian) Fock space
$F$ over a Hilbert space  $H$ is the completion of the tensor algebra
$$F=\oplus_{n=0}^{\infty}H^{\otimes n}.$$
Creation and annihilation operators are defined in the following way:
$$
A^{\dag}(f) f_{1}\otimes...\otimes f_{n}=f\otimes
f_{1}\otimes...\otimes f_{n}
$$
$$
A(f) f_{1}\otimes...\otimes
f_{n}=<f,f_{1}> f_{2}\otimes...\otimes f_{n}
$$
where  $<f,g>$ is the scalar product in the Hilbert space $H$.
Scalar product $\tau$ in the free Fock space is defined by the standard
construction of the direct sum of tensor products of Euclidean spaces.

We consider the case   $H=C\oplus C$, where $C$ is the field of
complex numbers.
In this case we have two creation operators
$A^{\dag}_{0}$, $A^{\dag}_{1}$
and two annihilation operators
$A_{0}$, $A_{1}$ with commutation relations
\begin{equation}\label{aac}
A_{i}A_{j}^{\dag}=\delta_{ij}.
\end{equation}
The vacuum vector $\Omega$  in the free Fock space satisfies
\begin{equation}\label{vacuum} A_{i}\Omega=0.
\end{equation}

We define
free coherent states   in the following way.

Let us fix two positive constants $\gamma_{0}$, $\gamma_{1}$,
$0<\gamma_{i}<1$.
Let us consider an infinite sequence of indices
$U=u_{1}u_{2}u_{3}...$, $u_i=0,1$.

The free coherent state $X_{U}$ is introduced as
the series
$$X_{U}=\sum_{k=0}^{\infty}X_{k}.$$
Here $X_{0}=\Omega$ is vacuum and
$X_{k}=\gamma_{u_{k}}A_{u_{k}}^{\dag}X_{k-1}$.
The norm  $||X_{k}||$ satisfies the condition
$$||X_{k}||=\prod_{i=1}^{k} \gamma_{u_{i}}$$
with $0<\gamma_{i}<1$.
Therefore the series $X_{U}$ converges.

Free coherent states are eigenvectors of the annihilation operator
$$
\gamma_{0}^{-1}A_{0}+\gamma_{1}^{-1}A_{1};
$$
for eigenvalue 1, i.e.
$$
(\gamma_{0}^{-1}A_{0}+\gamma_{1}^{-1}A_{1})X_{U}=X_{U}\qquad \forall U.
$$
Degeneration of this eigenspace is parametrized by
the set $\{U\}$ of infinite sequences of 0 and 1.

These sequences have  a natural interpretation as 2-adic numbers.
Every sequence $U$ is in one-to-one correspondence with a 2-adic number
$\sum_{i=1}^{\infty}u_i 2^{i-1}$.

Let us consider the following metric on the set of free
coherent states.
Let $U$, $V$ be arbitrary sequences of 0 and 1. These sequences
coincide up to an element with number $k$.
Let $X_{U}$, $X_{V}$ be corresponding free coherent states.
We introduce  the metric $\rho$:
$$
\rho(X_{U},X_{V})=||X_{k}||=\prod_{i=1}^{k} \gamma_{u_{i}}.
$$
The metric  $\rho$ is an ultrametric.
This means that the metric  $\rho$ obeys the strong triangle inequality
$$
\rho(X_{U},X_{V})\le \max(\rho(X_{U},X_{W}),\rho(X_{V},X_{W}))
$$
for arbitrary $U$, $V$, $W$.
Let us investigate the topology on the set of free coherent states
generated by the metric $\rho$.
Topology in a metric space is defined by the set of all closed balls
in this space.
The set of closed balls  $\{B_{U,k}\}$ in the metric space of free coherent
states with metric $\rho$ is parametrized by pairs $(U,k)$. Here $U$ is a
sequence of indices, $k$ is  non-negative integer and the ball $B_{U,k}$
contains all coherent states $X_{V}$ where $U$ and $V$ coincide up to the index
with number $k$.

The map $X_U\mapsto \sum_{i=1}^{\infty}u_i 2^{i-1}$ sends
the ball  $B_{U,k}$ to
a closed ball  of radius $2^{-k}$
with a center in $\sum_{i=1}^{\infty}u_i 2^{i-1}$
in the ring of integer 2-adic numbers.

Therefore the metric
$\rho(X_U,X_V)$
induces  the topology
of the ring of integer 2-adic numbers on the set of free coherent states.
The map $X_U\mapsto \sum_{i=1}^{\infty}u_{i}2^{i-1}$
is homeomorphism.

If
$\gamma_{0}=\gamma_{1}=\gamma$,
then $\rho(X_{U},X_{V})=\gamma^{k}$ is the metric on 2-adic disc.

Let us compare the metric  $\rho$ with the usual metric $\tau$
of the scalar product in the free Fock space. It is easy to see that
$$
\tau(X_{U},X_{V})=||X_{k}||\left(\sum_{i=k+1}^{\infty}
(\prod_{j=k+1}^{i}\gamma_{u_{j}}^{2}+
\prod_{j=k+1}^{i}\gamma_{v_{j}}^{2})\right)^{\frac{1}{2}}.
$$
Let  $\gamma_{0}<\gamma_{1}$. Then
$$
\left(\frac{2\gamma_{0}^2}{1-\gamma_{0}^2}\right)^{\frac{1}{2}}\rho
\le \tau \le
\left(\frac{2\gamma_{1}^2}{1-\gamma_{1}^2}\right)^{\frac{1}{2}}\rho;
$$
and metrics  $\tau$ and   $\rho$ generate equivalent topologies
on the space of free coherent states.
We have proved the following theorem.

{\bf   Theorem.}{\sl\qquad
The map $X_U\mapsto \sum_{i=1}^{\infty}u_{i}2^{i-1}$
is homeomorphism of the set of free coherent states with topology
generated by the metric of scalar product
$\tau(X_U,X_V)$
to  the ring of integer 2-adic numbers.
}

\vspace{3mm}
{\bf Acknowledgments}

Author is grateful to I.V.Volovich for discussions.

\end{document}